\documentclass[11pt,a4paper]{amsart}
\usepackage{graphicx,latexsym,amssymb}

\theoremstyle{plain}
\newtheorem{theorem*}{Theorem}
\newtheorem{theorem}{Theorem}[section]
\newtheorem{lemma}[theorem]{Lemma}

\newtheorem{proposition}[theorem]{Proposition}
\newtheorem{corollary}[theorem]{Corollary}
\newtheorem{corollary-lem}{Corollary}
\newtheorem{corollary*}[theorem*]{Corollary}
\theoremstyle{definition}
\newtheorem{definition}[theorem]{Definition}

\theoremstyle{remark}
\newtheorem{remark}[theorem]{Remark}

\numberwithin{equation}{section}
\numberwithin{figure}{section}
\font\nb=msbm10

\def \Z{\hbox{\nb  Z}}

\newcommand{\vb}{\overline{v}}
\newcommand{\wb}{\overline{w}}
\newcommand{\vbp}{\overline{v'}}
\newcommand{\wbp}{\overline{w'}}
\newcommand{\vbpp}{\overline{v''}}

\newcommand{\ben}{\begin{enumerate}}
\newcommand{\een}{\end{enumerate}}

\long\def\forget#1\forgotten{}
\begin{document}
\title{Palindromic Braids}
\date{\today}
\author[F. Deloup]{Florian Deloup}
\author[D. Garber]{David Garber}
\author[S. Kaplan]{Shmuel Kaplan}
\author[M. Teicher]{Mina Teicher}

\thanks{\tiny This paper is a part of the third author's Ph.D.\ Thesis at Bar-Ilan University.}
\thanks{\tiny The first author is a E.U. Marie
Curie Research Fellow (HMPF-CT-2001-01174). The second author is
partially supported by the Golda Meir Fellowship and wishes to
thank Ron Livne and the Einstein Institute of Mathematics in the
Hebrew University for hosting his stay. Third and Fourth authors are partially supported by EU-network HPRN-CT-2009-00099(EAGER),
Emmy Noether Research Institute for Mathematics, the Minerva
Foundation, and the Israel Science Foundation grant \#8008/02-3.}

\address{Florian Deloup,
Einstein Institute of Mathematics,
Edmond J. Safra Campus, Givat Ram,
The Hebrew University of Jerusalem,
91904 Jerusalem, Israel, {\it{and}}}
\address{Laboratoire Emile Picard,
UMR 5580 CNRS/Universit\'e Paul Sabatier, 118 route de Narbonne,
31062 Toulouse, France. } \email{deloup@picard.ups-tlse.fr}

\address{David Garber,
Einstein Institute of Mathematics,
Edmond J. Safra Campus, Givat Ram,
The Hebrew University of Jerusalem,
91904 Jerusalem, Israel, {\it{and}} }
\address{Department of Sciences,
Holon Academic Institute of Technology,
52 Golomb street,
58102 Holon, Israel.
}
\email{garber@math.huji.ac.il}
\email{garber@hait.ac.il}

\address{Shmuel Kaplan,
Department of Mathematics,
Bar-Ilan University
Ramat-Gan 52900,
Israel.
}
\email{kaplansh@macs.biu.ac.il}

\address{Mina Teicher,
Department of Mathematics,
Bar-Ilan University
Ramat-Gan 52900,
Israel.
}
\email{teicher@macs.biu.ac.il}

\begin{abstract}
The braid group $B_{n}$, endowed with Artin's presentation, admits
an antiautomorphism $B_{n} \to B_{n}$, such that $v \mapsto \vb$ is
defined by reading braids in reverse order (from right to left
instead of left to right). We prove that the map $B_{n} \to B_{n},
\ \ v \mapsto v \vb$ is injective. We also give some consequences arising due to this injectivity.
\end{abstract}
\maketitle 

\section{Introduction}

Let $n \geq 2$. Any free group $F_{n-1}$ on $n-1$ generators
$\sigma_1, \ldots, \sigma_{n-1}$ supports the antiautomorphism
$rev: w \mapsto \overline{w}$ defined by
$$ \sigma_{i_{1}}^{\alpha_{1}} \cdots \sigma_{i_{r}}^{\alpha_{r}}
\mapsto \sigma_{i_{r}}^{\alpha_{r}} \cdots
\sigma_{i_{1}}^{\alpha_{1}},$$
which reverses the
order of the word $w$ with respect to the prescribed set of
generators. It follows that any group $G$ presented by generators
and relations admits such an
antiautomorphism $rev$. The elements of $G$ which are
order-reversing invariant are called {\emph{palindromic}}.
In this
paper, we consider palindromic elements of Artin's Braid group
$B_{n}$, equipped with Artin's presentation, which will be called
{\emph{palindromic braids}}. Artin's presentation of the braid group $B_n$ consists of $n-1$
generators $\sigma_{1}, \ldots, \sigma_{n-1}$ and relations
\begin{equation}
\sigma_{i} \sigma_{j} = \sigma_{j} \sigma_{i} \ {\hbox{for}}\ \ \
|i - j | \geq 2, \label{eq:Artin1}
\end{equation}
\begin{equation}
\ \sigma_{i} \sigma_{i+1} \sigma_{i} = \sigma_{i+1} \sigma_{i}
\sigma_{i+1}\ \ \ {\hbox{for}}\ 1 \leq i \leq n-2.
\label{eq:Artin2}
\end{equation}

We distinguish between two equivalence relations on the elements of the braid group. For $a,b\in B_n$ we write $a=b$ to denote that $a$ and $b$ represent the same element in the group, and $a \equiv b$ to denote that $a$ and $b$ are actually the same element written letter by letter (i.e., $a\equiv b$ means that $a$ and $b$ are equal in the free group using only the generators of the braid group, with no relators).

\begin{figure}[htb]
\begin{center}
\includegraphics[height=5cm,width=8cm,angle=0,draft=false]{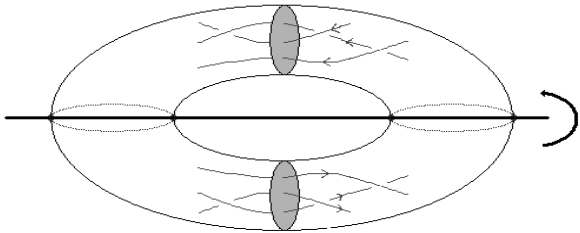}
\caption{The involution $inv$ and palindromic braids.}
\label{fig:hyperelliptic}
\end{center}
\end{figure}

Palindromic braids have a
particularly nice geometric interpretation. Given a geometric
braid $\beta$, denote by $\widehat{\beta}$ its closure into a link
inside a fixed solid torus $D^{2} \times S^{1}$. The solid torus
admits the involution
$$ {\rm{inv}}: D^{2} \times S^{1} \to D^{2} \times S^{1}, \ \ (re^{it},
\theta) \mapsto (r e^{-it}, - \theta),$$
whose set of
fixed points consists of two segments ($t \equiv 0$ (mod $\pi$) and $\theta
\equiv  0$ (mod $\pi$)), which is the intersection of the axis of the
180$^{\rm{o}}$ rotation with the solid torus. Observe that $\widehat{rev(\beta)}$ is
nothing else than inv$(\widehat{\beta})$ with the opposite
orientation. In particular, if a braid $\beta \in B_{n}$ is
palindromic then $\widehat{\beta}$ coincides with
${\rm{inv}}(\widehat{\beta})$ with the opposite orientation, see
Figure \ref{fig:hyperelliptic}.

We prove the following rigidity result for palindromic braids.

\begin{theorem} \label{th:main}
Let $\beta \in B_{n}$ be a palindromic braid such that
\begin{equation}
\beta = v \overline{v} \label{eq:decomposition-palindromic}
\end{equation}
for some braid $v \in B_{n}$. Then the decomposition
$(\ref{eq:decomposition-palindromic})$ is unique. Equivalently,
let $\beta = v \overline{v}$ and $\beta' = v' \overline{v'}$ be two
words in Artin's generators $\sigma_{1}, \ldots, \sigma_{n-1}$.
Then, $\beta = \beta'$ in $B_{n}$ if and only if $v = v'$ in $B_{n}$.
\end{theorem}

Of course one implication is obvious. Only the ``only if'' part of
the statement deserves a proof.

\begin{remark}
Note that Theorem \ref{th:main} cannot be generalized into the case of palindromic braid words of odd length. For example the two equal words $\sigma _1 \sigma _2 \sigma _1=\sigma _2 \sigma _1 \sigma _2$ are of the form $w \tau \overline{w}$ and $v \sigma \overline{v}$, however $\sigma _1=w \neq v=\sigma _2$. Moreover, not all palindromic braids of even length are of the form (\ref{eq:decomposition-palindromic}). For example, $\sigma _1 \sigma _3=\sigma _3 \sigma _1$ however, $\sigma _1 \neq \overline{\sigma _3}$.
\end{remark}

After this work was finished, F. Deloup communicated to us an alternative proof for Theorem \ref{th:main}, which is presented in \cite{Florian}, and is derived from the properties of the Dehornoy ordering of braids. The construction of the latter is a long process that requires rather sophisticated methods. In this paper we give an elementary proof, based on Garside normal form and its variant as developed by Jacquemard.

\section{Preliminaries and the Jacquemard Algorithm} \label{sec:generalities}

This section is devoted to the building blocks we use in order to prove
Theorem \ref{th:main}. Mainly, this section is intended to fix notations and recall some of the algorithms we use in this paper.

The monoid $B_{n}^{+}$ of \emph{positive braids} consists of braids which
admit a word representative which does not contain
$\sigma_{i}^{-1}$, $1 \leq i \leq n-1$.

Among positive braids, we can consider those
whose number of crossings between any two strands is less or equal to
$1$: they form the subset $S_{n}^{+} \subset B_{n}^{+}$ of {\emph{positive
permutation braids}}.

There is a canonical epimorphism $B_{n} \to S_{n}$. The image of a braid $\gamma$ is the permutation
{\emph{associated to}} $\gamma$. In particular, it is known that $S_{n}^{+}$ is in canonical bijection with the symmetric group $S_n$, which justifies the name of positive permutation braids.

There is only one positive braid $\Delta \in B_{n}^+$ in which any pair of strings crosses exactly once. It corresponds geometrically to a generalized half-twist which consists of all the strands $1,\cdots,n$, and is called the \emph{Garside element}. $\Delta$ is given by the formula:
\begin{equation}
\Delta =
(\sigma_{1}\sigma_{2} \ldots \sigma_{n-1})\
(\sigma_1\sigma_2\ldots \sigma_{n-2}) \cdots
(\sigma_{1}\sigma_{2})\ \sigma_1. \label{eq:Delta}
\end{equation}

A basic result asserts that the center of $B_{n}$ is generated by
$\Delta^2$. Abelianization of $B_{n}$ yields a canonical
homomorphism $B_{n} \to \Z$ which, when restricted to $B_{n}^{+}$,
coincides with word length with respect to Artin generators. We
denote by $|\beta|$ the length of $\beta \in B_n^+$; we have $|\sigma_{i}| = 1$ and the trivial braid
$e$ is the only positive braid such that $|e| = 0$.

We recall the algorithm given by Jacquemard \cite{EFF}, which manipulates a positive braid word $w \in B_n^+$ in order to write it using a given leading letter $\sigma _i$. The output of the algorithm is an equivalent positive braid word $\sigma _i w'=w$ or an indication that no $w'$ exists such that equality holds.

The basic nature of the algorithm is greedy. It starts by asking whether $w \equiv \sigma _i w'$ and stops if it does. If not, it looks for $\sigma _i$ inside $w$. In case $\sigma _i$ is not one of the letters of $w$, the algorithm returns \emph{false} which indicates non existence of $w' \in B_n^+$ such that $\sigma _i w'=w$.

When the algorithm found the leftmost $\sigma _i$ it works in two
steps: \ben \item Switch $\sigma _i$ with its left neighbor
$\sigma _j$ as long as $|i-j| \geq 2$. If $\sigma _i$ becomes the
first letter of the word we are done. However, in case that
$|i-j|=1$, the word is of the form $w=w_0\sigma _j \sigma _i w_1$,
and so in order to move $\sigma _i$ to the left one must use the
triple relation (\ref{eq:Artin2}) between $\sigma _j \sigma _i$
and the left most letter of $w_1$. If this is the case, the
algorithm goes to step $(2)$. \item The algorithm calls itself
recursively with the word $w_1$ and the letter $\sigma _j$. Upon
success of the recursive call the word looks like $w=w_0 \sigma _j
\sigma _i \sigma _j w_1'$ and therefore, we activate relation
(\ref{eq:Artin2}) on $\sigma _j \sigma _i \sigma _j$ resulting
with $w=w_0 \sigma _i \sigma _j \sigma _i w_1'$, and return to
step $(1)$. However, if the recursive call fails to extract
$\sigma _j$ to the left of $w_1$ the algorithm returns
\emph{false}. \een

To finish this section, we recall another result on the decomposition of braids, due to F. A. Garside \cite{Garside}
and later refined by W. P. Thurston \cite{Thurston}, and by E. A. Elrifai and H. R. Morton \cite{Elrifai-Morton}.

\begin{definition} \label{def:leftcan}
We say that a product $\alpha_1 \ldots \alpha_r$ satisfies {\emph{Thurston's condition}} if each $\alpha_{i}$ is a nontrivial positive permutation braid, and for any $1 \leq i \leq r-1$ we have that any $j$ such that $\alpha_{i+1}=\sigma _j \gamma _{i+1}$ also satisfies $\alpha _i=\gamma _i \sigma _j$ where $\gamma _i,\gamma _{i+1} \in S_n^+$.

\end{definition}

\begin{proposition}[Left-canonical form of a braid] \label{prop:leftcan}
Given any braid $\beta \in B_{n}$, there exists a unique decomposition
\begin{equation}
\beta = \Delta^{k} \alpha_{1} \cdots \alpha_{r}, \label{eq:can}
\end{equation}
where $k \in \Z$ is maximal, $\alpha_{i} \in S_{n}^{+}$ and the product
$\alpha_{1} \cdots \alpha_{r}$ satisfies Thurston's condition.
\end{proposition}

\section{Proof of Theorem \ref{th:main}}

Our first step consists of looking at the behavior of $\Delta$ and
permutation braids under the antiautomorphism $rev$.

\begin{lemma} \label{lem:basic}
The following properties hold:
\begin{enumerate}
\item[(1)] $\overline{\Delta} = \Delta$ \item[(2)] The set of
permutation braids is invariant under ${\rm{rev}}: w \mapsto
\overline{w}$. \item[(3)] $\overline{v^{-1}} = \overline{v}^{-1}$
for all $v \in B_{n}$.
\end{enumerate}
\end{lemma}

We now identify the basic problem. Let $\beta = \Delta^k \alpha_1 \ldots
\alpha_r$ be the left-canonical form for a braid $\beta \in B_{n}$.
We cannot assume that the decomposition $v = \alpha_{1} \ldots \alpha_{r}$ remains in
left-canonical form when viewed in $v \overline{v}$. Indeed, after multiplying on the right by the reversed braid,
the product $\alpha_{1} \ldots \alpha_{r}$ (viewed in
$v \overline{v}$) may cease to satisfy Thurston's condition.
 A simple example is provided by
$$ \alpha_1 = \sigma_1 \sigma_3, \ \ \alpha_2 = \sigma_3.$$
Both $\alpha_1$ and $\alpha_2$ are positive permutation braids and the product
$\gamma = \alpha_1 \alpha_2$ satisfies Thurston's condition. However, when we
write $\beta = \alpha _1 \alpha _2 \overline{\alpha}_2 \overline{\alpha}_1$ in
left-canonical form, we find $$\beta = \underbrace{\sigma_3
\sigma_1}_{\alpha_1}\
\underbrace{\sigma_3 \sigma_1}_{\alpha_2'}\ \underbrace{\sigma_3}_{\alpha_2''}
\underbrace{\sigma_3}_{\alpha_1''}$$
so that the second canonical factor $\alpha_2'$ does {\emph{not}}
coincide with $\alpha_{2}$.

We start by proving the Theorem \ref{th:main} for positive braids.

\subsection{Proof of Theorem \ref{th:main} for positive braid words}

We start by proving the following lemma:

\begin{lemma}\label{lem:positive}
Let $v \vb=\sigma w' \wbp \sigma$ and let $v \vb\equiv w_0=w_1=\cdots =w_k\equiv \sigma v'$ be a sequence of positive braid words such that each $w_{i+1}$ is the outcome of the activation of one relation out of the relations in the semigroup $B_n^+$ on $w_i$ according to Jacquemard's algorithm . Then, all relations are performed only within the first half of the word $w_i$ which implies they all involve only letters from $v$.
\end{lemma}

\begin{proof}
Notice that since $v \vb=\sigma w' \wbp \sigma$, the success of the Jacquemard's algorithm is guaranteed. Hence we know that $\sigma $ is one of the letters of $v$. Now, we need to prove that in each step of Jacquemard's algorithm that uses a relation, it occurs in the first half of the word $v$.

For step $(1)$ of the algorithm, this is obvious: all relations involve $\sigma$ and \emph{left} neighbors of $\sigma$; since $\sigma$ is in $v$, all relations occur inside $v$. Moreover, relations can be activated mirror like on $\vb$ as well; Hence, we maintain the palindromic structure of the word. This implies that when we need to move to step $(2)$ of the algorithm we have
$$v\vb=v_1 \tau \sigma v_2 \vb = v_1 \tau \sigma v_2 \overline{v_2} \sigma \tau \overline{v_1},$$
where $\tau$ and $\sigma$ do not commute, and $\sigma$ is not in $v_1$.

If this is the case, the algorithm calls itself recursively using $v_2 \vb=v_2 \overline{v_2}\sigma \tau \overline {v_1}$ and $\tau$, trying to extract $\tau$ to the left, first by looking for the leftmost $\tau$ letter in $v_2 \vb=v_2 \overline{v_2}\sigma \tau \overline {v_1}$.

Assume by contradiction that $\tau$ is not a letter of $v_2$ (hence is not a letter of $\overline{v_2}$). Then the leftmost $\tau$ letter in $v_2 \vb=v_2 \overline{v_2}\sigma \tau \overline {v_1}$ appears to the right of $\sigma$ and to the left of $\overline{v_1}$. In order to extract this $\tau$ to the left of $v_2$, we need to activate another recursive call of the algorithm on $\overline{v_1}$ with the letter $\sigma$ (since in our case $\sigma$ and $\tau$ do not commute and we have to use step $(2)$ of the algorithm). But, $\sigma$ is not in $\overline{v_1}$ (since it is not in $v_1$). Therefore $\sigma$ cannot be extracted to the left of $\overline{v_1}$. This implies that Jacquemard's algorithm failed, and this is a contradiction.

If $\tau$ is found within $v_2$ and during the process of extracting it to the left we do not encounter the need to use letters from the right half of $v\vb$ we are finished. Therefore, assume by contradiction that at some point in the process we encounter a relation involving a letter of $\vb$. Again, since until this step all relations were activated only in the left half of the word their mirror image can be activated on the right half of the word, so the palindromic structure of the word is preserved.

Suppose that we performed $k$ recursive steps of the algorithm. Then, our word looks like:

$$v\vb=v_1\tau _1\sigma v_2 \tau _2 \tau _1 v_3 \tau _3 \tau _2 \cdots v_{k-1} \tau _{k-1} \tau _{k-2} v_k \tau _k \tau _{k-1} v_{k+1}\ \cdot$$
$$\overline{v_{k+1}} \tau _{k-1} \tau _k \overline{v_k} \tau _{k-2} \tau _{k-1} v_{k-1} \cdots \tau _2 \tau _3 \overline{v_3} \tau _1 \tau _2 \overline{v_2} \sigma \tau _1 \overline{v_1}$$
where $\tau _i$ and $\tau _{i+1}$ are two non commuting letters, $\sigma$ does not commute with $\tau _1$ and is not a letter of $v_1$. Moreover, $\tau _i$ is not a letter of $v_{i+1}$ for any $i=1,\cdots, k$.

Now, in this recursion step, we have called the algorithm with the letter $\tau _k$ and the word $v_{k+1} \overline{v_{k+1}} \tau _{k-1} \tau _k \overline{v_k} \tau _{k-2} \tau _{k-1} v_{k-1} \cdots \tau _2 \tau _3 \overline{v_3} \tau _1 \tau _2 \overline{v_2} \sigma \tau _1 \overline{v_1}$. However, since $\tau _k$ is not a letter of $v_{k+1}$ and of $\overline{v_{k+1}}$ the leftmost $\tau _k$ in this recursion call is to the right of $\tau _{k-1}$ and to the left of  $\overline{v_k}$. Since $\tau _{k-1}$ does not commute with $\tau _k$, another recursion call is needed with the letter $\tau _{k-1}$ and the word $\overline{v_k} \tau _{k-2} \tau _{k-1} v_{k-1} \cdots \tau _2 \tau _3 \overline{v_3} \tau _1 \tau _2 \overline{v_2} \sigma \tau _1 \overline{v_1}$. Again, $\tau _{k-1}$ is not a letter of $\overline{v_k}$ hence the leftmost $\tau _{k-1}$ in this recursion call is to the right of $\tau _{k-2}$ and to the left of $\overline{v_{k-1}}$. Similarly $\tau _{k-2}$ does not commute with $\tau _{k-1}$, so we continue $k-3$ recursion calls until we reach a recursion call with the letter $\tau _1$ and the word $\overline{v_2} \sigma \tau _1 \overline{v_1}$. Since $\tau _1$ is not a letter in $\overline{v_2}$, the leftmost $\tau _1$ in this call is to the right of $\sigma$ and to the left of $\overline{v_1}$. This implies that another recursion call is needed in order to extract the letter $\sigma$ from the word $\overline{v_1}$. However, this contradicts the hypothesis on $v_1$.

This concludes the proof of all cases, hence all relations are activated inside the left half of the word $v \vb$ as claimed.
\end{proof}

Now we are ready to prove the theorem for positive braid words.

\begin{theorem} \label{th:positive}
Let $\beta , \beta ' \in B_{n}^+$ be two palindromic positive braids of even length such that
$\beta = v \vb$ and $\beta '= w \wb$ for some braids $v,w \in B_{n}^+$,
Then, $\beta = \beta'$ in $B_{n}$ if and only if $v = w$ in $B_{n}$.
\end{theorem}

\begin{proof}

By induction on the length $l=|w|$ of $w$.
Assume that $w \equiv \sigma w'$ where $w' \in B_n^+$, i.e., $\sigma$ is the first letter in $w$. Then, $w \wb=\sigma w' \wbp \sigma$. This means that $v \vb$ can be written such that its first letter is $\sigma$, that is, $v \vb=\sigma v'$ for some $v' \in B_n^+$.

By the embedding theorem of Garside \cite{Garside} it follows that there is a sequence of words $v \vb\equiv w_0=w_1=\cdots =w_k\equiv \sigma v'$, such that each $w_{i+1}$ is obtained from $w_i$ by activating one relation out of the relations in the semigroup $B_n^+$.

One possible sequence is the one which uses the relations suggested by the algorithm of Jacquemard given in \cite{EFF}. Now, Lemma \ref{lem:positive} shows that every relation used in the sequence is fully contained in $v$ (the left half of the word), and does not effect $\vb$. Therefore, it is possible to activate all relations described in the sequence in a mirror-like image on $\vb$ and get $w\wb=\sigma w'\wbp \sigma = v \vb = \sigma v' = \sigma v'' \vbpp \sigma$.

However, if this is the case $w'\wbp=v''\vbpp$ where $|w'\wbp|=|v''\vbpp|=l-2$. Therefore, by the induction hypothesis we have $w'=v''$ which implies that $\sigma w'=w=v=\sigma v''$.

\end{proof}

Next, we use the above, and give the proof of Theorem \ref{th:main}.

\subsection{Proof of Theorem \ref{th:main} for the general case}

Let $v \in B_{n}$ and set $\beta = v \vb$.
We need to prove that if $\beta = w \wb$, then $v=w$.

Let $v = \Delta^{k} \alpha_{1} \ldots \alpha_{r}$ be the left-canonical form of $v$.
Then  $\overline{v} = \overline{\alpha}_{r} \ldots \overline{\alpha}_{2} \overline{\alpha}_{1}\ \Delta^{k}$ (Lemma \ref{lem:basic}$(1)$). Hence,
 $$
 \beta = v \overline{v} = \Delta^{k} \alpha_{1} \ldots \alpha_{r}
 \overline{\alpha}_{r} \ldots
 \overline{\alpha}_{1}  \Delta^{k}.
 $$

Moreover, let $w=\Delta ^j \beta _1 \cdots \beta _p$ be the left-canonical form of $w$. We have

$$\beta = w\wb=\Delta ^j \beta _1 \cdots \beta _p \overline{\beta} _p \cdots \overline{\beta} _1 \Delta ^j$$
(Note that it is not necessary that $r=p$).

Without loss of generality, we may assume that $j<k<0$ (Otherwise, $\beta \in B_n^+$ and we use Theorem \ref{th:positive}). By multiplying $\beta$ by $\Delta ^{-j}$ on the left and on the right, we obtain:

$$\Delta ^{-j} \beta \Delta ^{-j}=\Delta ^{-j} v\vb \Delta ^{-j}=\Delta ^{k-j} \alpha _1 \cdots \alpha _r \overline{\alpha} _r \cdots \overline{\alpha _1} \Delta ^{k-j},$$
and
$$\Delta ^{-j} \beta \Delta ^{-j}=\Delta ^{-j} x\overline{x} \Delta ^{-j}=\beta _1 \cdots \beta _p \overline{\beta} _p \cdots \overline{\beta _1}.$$
These are two equal positive braid words and have the form $v'\vbp=w'\wbp$, where $v'=\Delta ^{k-j} \alpha _1 \cdots \alpha _r$ and $w'=\beta _1 \cdots \beta _p$ respectively. Therefore, Theorem \ref{th:positive} applies and we conclude that $v=w$.

\hfill
$\qed$

As a consequence, we obtain the following corollaries:

\begin{corollary} \label{co:1stcorollary}
Let $\beta =x\overline{x} \in B_n^+$ be a positive palindromic braid of even length, and let the left-canonical normal form of $x$ be $\alpha _1 \cdots \alpha _r$ such that $\alpha _1 \neq \Delta$. Then, the left-canonical normal form of $x\overline{x}$ is $\beta _1 \cdots \beta _p$ where $\beta _1 \neq \Delta$.
\end{corollary}

\begin{proof}
Otherwise because the process of Jacquemard's algorithm may be used to transform the word $\alpha_1 \cdots \alpha_r\overline{\alpha}_r\cdots\overline{\alpha}_1$ into its left-canonical form, and since it extracts letters to the left only from the first half of the word, we might have
$$\beta=\alpha _1 \cdots \alpha _r \overline{\alpha} _r \cdots \overline{\alpha} _1=\Delta \gamma _1 \cdots \gamma _q \overline{\gamma} _q \cdots \overline{\gamma} _1 \Delta,$$
which means, by Theorem \ref{th:positive}, that $\alpha _1 \cdots \alpha _r = \Delta \gamma _1 \cdots \gamma _q$. By the uniqueness of the left-canonical normal form, we deduce that $\alpha _1=\Delta$, which is a contradiction.
\end{proof}

We may generalize Corollary \ref{co:1stcorollary}:
\begin{corollary}
Let $v \in B_n^+$ be a positive braid, and let $n(v)$ denote the number of leading permutation braids which are $\Delta$ when $v$ is written in its left-canonical normal form. Then, for $\beta=v \vb$, we have $n(v\vb)=2n(v)$.
\end{corollary}

\begin{proof}
Since $v$ is written in left-canonical form as $\Delta ^{n(v)} \alpha _1 \cdots \alpha _r$ where $\alpha _1 \neq \Delta$, we have that $\beta=\Delta ^{n(v)} \alpha _1 \cdots \alpha _r \overline{\alpha}_r \cdots \overline{\alpha}_1 \Delta ^{n(v)}$. Note that since $\Delta$ almost commutes with any permutation braid, we may write
$$\beta=\Delta ^{2n(v)} \alpha _1' \cdots \alpha _r' \overline{\alpha}'_r \cdots \overline{\alpha}'_1,$$ where $\alpha _i'=\alpha _i$ if $n(v)$ is even and $\alpha _i'$ is obtained from $\alpha _i$ by replacing each $\sigma _j$ by $\sigma _{n-j}$ in case that $n(v)$ is odd. In any of the cases, the product $\alpha _1' \cdots \alpha _r'$ keeps it's left-canonical form. Hence, using the same argument as in Corollary \ref{co:1stcorollary} we get that $n(v\overline{v})=2n(v)$.
\end{proof}

\bibliographystyle{amsalpha}

%
%
%
%
%

\vskip 2cm

\end{document}